\documentclass[12pt]{amsart}
\usepackage{amsmath,amssymb,amscd,latexsym}

\newtheorem{Theorem}{Theorem}[section]

\newtheorem{Proposition}[Theorem]{Proposition}
\newtheorem{Lemma}[Theorem]{Lemma}

\newtheorem{Remark}[Theorem]{Remark}

\begin{document}

\title{On Backward Uniqueness\\ for the Heat Operator in Cones}

\author{Jie Wu \,and\, Wendong Wang}

\begin{figure}[b]

{\small
\begin{tabular}{ll}
\end{tabular}}
\end{figure}

\date{June, 2013}

\begin{abstract}
Consider the system $|\partial_tu+\Delta u|\leq M(|u|+|\nabla u|)$,
$|u(x,t)|\leq Me^{M|x|^2}$ in $\mathcal{C}_{\theta}\times[0,T]$ and $u(x,0)=0$ in $\mathcal{C}_{\theta}$, where $\mathcal{C}_{\theta}$ is a cone with opening angle $\theta$. L. Escauriaza constructed an example to show that such system has a nonzero bounded solution when $\theta<90^\circ$, and it's conjectured that the system has only zero solution for $\theta>90^\circ$. Recently Lu Li and V. \v{S}ver\'{a}k \cite{LlS} proved that the claim is true for $\theta>109.5^\circ$. Here we improve their result and prove that only zero solution exists for this system when $\theta>99^\circ$ by exploring a new type of  Carleman inequality, which is of independent interest.
\end{abstract}

\maketitle

\\


\section{Introduction}
Let $U$ be a domain in $\mathbb{R}^n$ and $u$ be the solution of the following equation:
\begin{equation}\label{classical1}
\partial_tu-\Delta u+b(x,t)\cdot\nabla u+c(x,t)u=0, ~~~(x,t)\in U\times(0,T),
\end{equation}
where $b$ and $c$ are bounded.  Moreover, assume that $u$ satisfies the growth condition
$$|u(x,t)|\leq Me^{M|x|^2},$$
for some $M>0$. The backward uniqueness (BU) problem is: if $u$ vanishes at $t=T$, that is
$$u(x,T)=0,~~~x\in U,$$
does $u$ vanish identically in $U\times(0,T)$? If so, we say that $U$ is a (BU) domain. The point is that no any parabolic boundary conditions are made upon  $u$.

The backward uniqueness property is applied in many problems, for example, the control theory for PDEs and the
regularity theory of parabolic equations. Especially, it plays an important role in the proof of critical $L^{\infty,3}_{t,x}$ regularity for 3D Navier-Stokes equations, see \cite{ESS2}.

C. C. Poon proved that the whole space is a (BU) domain in \cite{Po}. L. Escauriaza, G. Seregin and V. \v{S}ver\'{a}k proved that the exterior of a ball \cite{ESS} and the half space \cite{ESS2,ESS3} are both (BU) domains. The tools they used to prove (BU) are frequency function methods and Carleman inequality methods. On the other hand, any bounded domain is not (BU) domain, see \cite{Jone,Lit}. For parabolic operators of variable, singular or degenerated coefficients on the whole space, we refer to \cite{LPR,Ho,LA,LIN,Ch,SP,SP2,WZ,WUZ} and references therein for more relevant interesting results.

In this paper, our attention is focused on the (BU) problem in cones, which seems to be rather interesting. In \cite{LlS}, L. Escausriaza constructed an example to show that (BU) fails when the opening angle $\theta$ of the cone satisfies $\theta<\frac{\pi}{2}$, and $\theta=\frac{\pi}{2}$ seems to be the borderline case for Escausriaza's construction. Recently, Lu Li and V. \v{S}ver\'{a}k proved that (BU) holds when $109.5^\circ<\theta<\pi$. Inspired by the above results, it is conjectured naturally that

\emph{ (BU) holds when $\frac{\pi}{2}<\theta<\pi$ and fails when $\theta<\frac{\pi}{2}$. }

 Here we improve the result of Lu Li and V. \v{S}ver\'{a}k and obtain that (BU) holds when $98.99^\circ<\theta<\pi$.

We consider the following system
\begin{eqnarray}\label{system}
\left\{
  \begin{array}{ll}
    |\partial_tu+\Delta u|\leq M(|u|+|\nabla u|), \qquad &{\rm in} \quad \mathcal{C}_{\theta}\times(0,T),  \\
    |u(x,t)|\leq Me^{M|x|^2}, \qquad &{\rm in} \quad \mathcal{C}_{\theta}\times(0,T),  \\
    u(x,0)=0, \qquad &{\rm in} \quad \mathcal{C}_{\theta},
  \end{array}
\right.
\end{eqnarray}
where $\mathcal{C}_{\theta}$ is a cone with opening angle $\theta$. In a suitable coordinates, we can write $\mathcal{C}_{\theta}$ as
$$\mathcal{C}_{\theta}=\{x=(x_1,x')|x'\in\mathbb{R}^{n-1},x_1>|x|\cos{\frac{\theta}{2}}\}.$$

Let $\varepsilon=\cos{\frac{\theta}{2}}<\varepsilon_0$ and our goal is the border line $\varepsilon_0\rightarrow\sqrt{\frac12}\approx0.7070$. In fact, we obtain that  $\varepsilon_0\approx 0.6495$, $2\arccos{\varepsilon_0}\approx98.99^\circ$ and $m\approx2.4600$, where $\varepsilon_0,m$ are chosen the biggest numbers such that for $\frac12<\gamma\leq 1$, $2.36\leq m<3$, there hold (see the estimate of $J_3$ in Section 3)
\begin{eqnarray}\label{eq:j3 condition}
\left\{\begin{array}{lll}
&(i)~4(2\gamma-1)=\gamma^2(4-\gamma^2\frac{m-1}{4});\\
&(ii)~\frac{m-1}{m+1}-\varepsilon_0^2= 0;\\
&(iii)~(4-\gamma^2\frac{m-1}{4}-m)-(4-\gamma^2\frac{m-1}{4})\varepsilon_0^m=0.
\end{array}\right.
\end{eqnarray}

Our main result is the following
\begin{Theorem}\label{thm}
Assume that $u$ satisfies (\ref{system}) for some $\theta\in(2\arccos{\varepsilon_0},\pi)$ with $\varepsilon_0\approx 0.6495$,
then $u(x,t)\equiv0$ in $\mathcal{C}_{\theta}\times(0,T)$.
\end{Theorem}

The proof of the above theorem is based on the following Carleman inequality. Assuming it, there is only a standard argument left to prove Theorem \ref{thm}.
\begin{Proposition}\label{Prop}
For $m$ and $\varepsilon_0$ as above, set $\varepsilon\in(0,\varepsilon_0)$ and $\varepsilon=\cos{\frac{\theta}{2}}$. Moreover, we assume that
$$Q=\{(x,t)|x\in \mathcal{C}_{\theta},x_1>1,t\in(0,1)\},\quad
\varphi(x)=r^\alpha[(\frac{x_1}{r})^m-\varepsilon^m]$$
where
$r=|x|$. Then there exist $\alpha=\alpha(\varepsilon)\in(1,2)$ and $K=K(\varepsilon)>0$, such that  the following inequality holds for any $u\in C_0^\infty(Q)$ and any real number $a>0$:
\begin{equation}\label{CI}
\begin{split}
\int_Qe^{2a(t^{-K}-1)\varphi(x)-\frac{|x|^2+K}{8t}}(|u|^2+|\nabla u|^2)dxdt
\leq \int_Qe^{2a(t^{-K}-1)\varphi(x)-\frac{|x|^2+K}{8t}}|\partial_tu+\Delta u|^2dxdt.
\end{split}
\end{equation}
\end{Proposition}

\begin{Remark}
Theorem \ref{thm} improves previous results obtained in \cite{LlS}. Furthermore, a new  Carleman inequality  (\ref{CI}) is proved for Theorem \ref{thm},
which  is somewhat different from that given by Lu Li and V. \v{S}ver\'{a}k in \cite{LlS}.
In the past, we need two Carleman inequalities to prove (BU) (for example, see \cite{ESS3,LlS}): the first one appears with the heat kernel weight, which implies an
exponential decay of the solution with respect to time at $x=0$, and the second owns the weight about the cone boundary property. However, here the Carleman inequality (\ref{CI}) whose weight owns the above two properties is sufficient
to prove (BU).
\end{Remark}

The paper is organized as follows. First we prove Theorem \ref{thm} in Section 2 under the assumption of Carleman inequality (\ref{CI}), then the next section is devoted to the proof of Carleman inequality (\ref{CI}).

\section{Proof of Theorem \ref{thm}}
Without loss of generality, assume that $T=1$. We always extend $u(x,t)$ by zero to the negative values of $t$.\\

For $x\in\mathcal{C}_{\theta}$, we denote the distance between $x$ and the boundary of $\mathcal{C}_{\theta}$ by $d_\theta(x)$ as in \cite{LlS}, then
$$d_\theta(x)=x_1\sin\frac{\theta}{2}-|x'|\cos\frac{\theta}{2}.$$
For $b>0$, we denote
$$\mathcal{C}_{\theta}^{+b}=\{x\in\mathcal{C}_{\theta}| \,d_\theta(x)>b\}.$$
The proof of Theorem \ref{thm} is based on the following lemma.
\begin{Lemma}\label{lem}
Assume that $u$ satisfies (\ref{system}) for some $\theta\in(2\arccos{\varepsilon_0},\pi)$,
then there exists $T_1=T_1(M)\in(0,\frac{1}{2})$, such that
$$u(x,t)\equiv0$$
in $\mathcal{C}_{\theta}\times(0,T_1)$.
\end{Lemma}
\emph{Proof}.
We'll use Carleman inequality (\ref{CI}) to prove this lemma, and mainly follow the same line as in \cite{ESS,ESS3,LlS}.

At first, by the regularity theory for solutions of parabolic equations \cite{LS}, we have
\begin{equation}\label{pt1.1}
|u(x,t)|+|\nabla u(x,t)|\leq C(M)e^{2M|x|^2}
\end{equation}
for all $(x,t)\in \mathcal{C}_{\theta}^{+1}\times(0,\frac{1}{2})$.

Let
\begin{equation}\label{T_1}
T_1=\min\{\frac{1}{256M},\frac{1}{12M^2},\frac{1}{2}\},
\end{equation}
$\lambda=\sqrt{2T_1}$, and $y_0=(1/\sin\frac{\theta}{2},0,\cdots, 0)$.  For $(y,s)\in \mathcal{C}_\theta\times(0,1)$, define that
$$v(y,s)=u(\lambda y+y_0,\lambda^2s-T_1).$$
By (\ref{system}) we have
\begin{equation}\label{pt1.2}
|\partial_sv+\Delta v|\leq\lambda M(|v|+|\nabla v|),
\end{equation}
for all $(y,s)\in\mathcal{C}_\theta\times(0,1)$.

Note that for $(y,s)\in\mathcal{C}_\theta\times(\frac12,1)$, we get $$(\lambda y+y_0,\lambda^2s-T_1)\in\mathcal{C}_\theta^{+1}\times(0,T_1),$$ then (\ref{pt1.1}) implies that
\begin{equation}\label{pt1.3}
|v(y,s)|+|\nabla v(y,s)|\leq C(M)e^{2M|\lambda y+y_0|^2}
\leq C(M,\varepsilon)e^{4M\lambda^2|y|^2},
\end{equation}
and
\begin{equation}\label{pt1.4}
v(y,s)=0,\quad {\rm for }~~(y,s)\in \mathcal{C}_\theta\times(0,\frac{1}{2}].
\end{equation}

In order to apply Carleman inequality (\ref{CI}), we choose two smooth cut-off functions $\eta_1,\eta_2$ satisfying
$$
\eta_1(p)=\left\{
                 \begin{array}{ll}
                   0, & \hbox{if $p<2$,} \\
                   1, & \hbox{if $p>3$,}
                 \end{array}
               \right.
$$
and
$$
\eta_2(q)=\left\{
            \begin{array}{ll}
              0, & \hbox{if $q<-\frac{3}{4}$,} \\
              1, & \hbox{if $q>-\frac{1}{2}$.}
            \end{array}
          \right.
$$
Furthermore, $0\leq \eta_1,\eta_2\leq 1$; $|\eta_1'|$, $|\eta_1''|$, $|\eta'_2|$ and $|\eta''_2|$ are all bounded. Denote $\Lambda(s)=s^{-K}-1$ and
$$C_\star=1+\sup_{\frac{1}{2}<s<1 \atop y\in\mathcal{C}_\theta,2<y_1<3}\{\Lambda(s)\varphi(y)\},$$
where $\varphi(y)$ is defined in Proposition 1.2, and $C_\star$ is well-defined since $m>\alpha.$

Let $$\eta(y,s)=\eta_1(y_1)\eta_2(\frac{\Lambda(s)\varphi(y)}{2C_\star}-1),$$
and $w=\eta v$, then $supp~w\subset Q$ and
\begin{eqnarray*}
\begin{split}
|\partial_sw+\Delta w|=&|\eta(\partial_sv+\Delta v)+v(\partial_s\eta+\Delta \eta)+2\nabla v\nabla\eta|\\
\leq&|\eta(\partial_sv+\Delta v)|+2\chi_{\Omega}(|v|+|\nabla v|)(|\partial_s\eta|+|\nabla\eta|+|\nabla^2\eta|),
\end{split}
\end{eqnarray*}
where $\chi$ is the characteristic function of the set
$$\Omega=\{(y,s)|\frac{1}{2}<s<1,0<\eta(y,s)<1\}.$$
By (\ref{pt1.2}) we have
\begin{eqnarray*}
\begin{split}
|\partial_sw+\Delta w|\leq& \eta\lambda M(|v|+|\nabla v|)+2\chi(|v|+|\nabla v|)(|\partial_s\eta|+|\nabla\eta|+|\nabla^2\eta|)\\
\leq& \lambda M(|w|+|\nabla w|)+C(M)\chi_{\Omega}(|v|+|\nabla v|)(|\partial_s\eta|+|\nabla\eta|+|\nabla^2\eta|).\\
\end{split}
\end{eqnarray*}
Note that $\frac{1}{2}<s<1$ in $\Omega$ and $|\nabla^k\varphi(y)|\leq C(\varepsilon)|y|^{\alpha-k}$ with $ k=1,2$, then when $\frac{1}{2}<s<1$,
$$(|\partial_s\eta|+|\nabla\eta|+|\nabla^2\eta|)\leq C(\varepsilon)|y|^{\alpha}\leq C(\varepsilon)|y|^{2}.$$
Thus
\begin{equation}\label{pt1.5}
|\partial_sw+\Delta w|\leq \lambda M(|w|+|\nabla w|)+C(M,\varepsilon)\chi_{\Omega}(|v|+|\nabla v|)|y|^2.
\end{equation}
Especially,
\begin{eqnarray*}
\begin{split}
\Omega=&\{(y,s)|\frac{1}{2}<s<1,\eta_1>0,0<\eta_2<1\}\\
&\bigcup\{(y,s)|\frac{1}{2}<s<1,0<\eta_1<1,\eta_2=1\}\\
=&\{(y,s)|\frac{1}{2}<s<1,y_1>2,\frac{1}{2}<\frac{\Lambda(s)\varphi(y)}{C_\star}<1\}\\
&\bigcup\{(y,s)|\frac{1}{2}<s<1,2<y_1<3,\frac{\Lambda(s)\varphi(y)}{C_\star}\geq1\},
\end{split}
\end{eqnarray*}
Due to the choice of $C_\star$, we obtain that the second set of the right-hand side in the above identity is empty, then
\begin{equation}\label{pt1.6}
\Omega=\{(y,s)|\frac{1}{2}<s<1,y_1>2,\frac{1}{2}<\frac{\Lambda(s)\varphi(y)}{C_\star}<1\}.
\end{equation}
By (\ref{pt1.3}), we derived that
\begin{eqnarray*}
\begin{split}
e^{2a \Lambda(s)\varphi(y)-\frac{|y|^2+K}{8s}}(|v|+|\nabla v|)^2
\leq& C(M,\varepsilon)e^{2a \Lambda(s)\varphi(y)-\frac{|y|^2+K}{8s}+8M\lambda^2|y|^2}\\
\leq& C(M,\varepsilon)e^{2a \Lambda(s)|y|^\alpha-\frac{|y|^2}{8}+16MT_1|y|^2}\\
\end{split}
\end{eqnarray*}
With the help of (\ref{T_1}) and $T_1\leq\frac{1}{256M}$, we have in $Q$
\begin{equation}\label{pt1.7}
e^{2a \Lambda(s)\varphi(y)-\frac{|y|^2+K}{8s}}(|v|+|\nabla v|)^2\leq C(M,\varepsilon)e^{2a \Lambda(s)|y|^\alpha-\frac{|y|^2}{16}}.
\end{equation}
Although $supp~w$ may be unbounded, $supp~w\subset Q$ and (\ref{pt1.7}) allow us to claim the validity of Proposition \ref{Prop} for $w$. Then by  Carleman inequality (\ref{CI}), together with (\ref{pt1.5}), we get
\begin{eqnarray*}
\begin{split}
J&\equiv\int_Qe^{2a \Lambda(s)\varphi(y)-\frac{|y|^2+K}{8s}}(|w|^2+|\nabla w|^2)dyds\\
&\leq \int_Qe^{2a \Lambda(s)\varphi(y)-\frac{|y|^2+K}{8s}}|\partial_sw+\Delta w|^2dyds\\
&\leq 3\lambda^2M^2J+C(M,\varepsilon)\int_Qe^{2a \Lambda(s)\varphi(y)-\frac{|y|^2+K}{8s}}\chi_{\Omega}(|v|+|\nabla v|)^2|y|^4dyds.
\end{split}
\end{eqnarray*}
Using (\ref{T_1}) again, we know that $3\lambda^2M^2=6T_1 M^2\leq\frac{1}{2}$. The above inequality implies that
$$J\leq C(M,\varepsilon)\int_\Omega e^{2a \Lambda(s)\varphi(y)-\frac{|y|^2+K}{8s}}(|v|+|\nabla v|)^2|y|^4dyds,$$
and by (\ref{pt1.7}) we have
$$J\leq C(M,\varepsilon)\int_\Omega e^{2a \Lambda(s)\varphi(y)-\frac{|y|^2}{16}}|y|^4dyds.$$
Since $\Lambda(s)\varphi(y)<C_\star$ in $\Omega$ according to (\ref{pt1.6}) , then
\begin{equation}\label{pt1.8}
J\leq C(M,\varepsilon) e^{2a C_\star}\int_\Omega e^{-\frac{|y|^2}{16}}|y|^4dyds
\leq C(M,\varepsilon) e^{2a C_\star}.
\end{equation}

On the other hand, note that $\varphi(y)\rightarrow\infty$ as $|y|\rightarrow\infty$, then we define
$$\Omega_1=\{(y,s)|0<s<1,\eta=1\}=\{(y,s)|0<s<1, y_1\geq3, \frac{\Lambda(s)\varphi(y)}{C_\star}\geq1\},$$
and
\begin{equation}\label{omega2}
\Omega_2=\{(y,s)|0<s<1, y_1\geq3, \frac{\Lambda(s)\varphi(y)}{C_\star}\geq2\}\neq \emptyset.
\end{equation}
Obviously $\Omega_2\subset\Omega_1$ and $w=v$ in $\Omega_1$. Hence
\begin{eqnarray*}
\begin{split}
J&\geq\int_{\Omega_1}e^{2a \Lambda(s)\varphi(y)-\frac{|y|^2+K}{8s}}(|v|^2+|\nabla v|^2)dyds\\
&\geq\int_{\Omega_2}e^{2a \Lambda(s)\varphi(y)-\frac{|y|^2+K}{8s}}(|v|^2+|\nabla v|^2)dyds.
\end{split}
\end{eqnarray*}
From (\ref{omega2}) we know that $\Lambda(s)\varphi(y)\geq2C_\star$ in $\Omega_2$, hence
\begin{equation}\label{pt1.9}
J\geq e^{4a C_\star}\int_{\Omega_2}e^{-\frac{|y|^2+K}{8s}}(|v|^2+|\nabla v|^2)dyds.
\end{equation}

Combining (\ref{pt1.8}) and (\ref{pt1.9}), finally we obtained that
$$\int_{\Omega_2}e^{-\frac{|y|^2+K}{8s}}(|v|^2+|\nabla v|^2)dyds\leq C(M,\varepsilon)e^{-2a C_\star}.$$
Passing to the limit as $a\rightarrow+\infty$, we obtain $v(y,s)=0$ in $\Omega_2$. Using unique continuation though spatial
boundaries (see \cite{ESS2}), we obtain that $v(y,s)=0$ in $\mathcal{C}_\theta\times(0,1)$, then $u(x,t)=0$ in $\mathcal{C}_\theta^{+2}\times(0,T_1)$. Using the unique continuation result again , we get $u(x,t)=0$ in $\mathcal{C}_\theta\times(0,T_1)$. Thus we have proved this lemma.\\

Nest we give the complete proof of Theorem \ref{thm}.\\

\emph{Proof of Theorem \ref{thm}}. We define
$$u^{(1)}(y,s)=u(\sqrt{1-T_1}y,(1-T_1)s+T_1),~~~(y,s)\in\mathcal{C}_\theta\times(0,1).$$
Then $u^{(1)}$ satisfies the conditions of Lemma \ref{lem}, and we have
$u^{(1)}(y,s)=0$ in $\mathcal{C}_\theta\times(0,T_1)$.\\
In other words, $u(x,t)=0$ in $\mathcal{C}_\theta\times(0,T_2)$, where
$$T_2=(1-T_1)T_1+T_1.$$
After iterating $k$ steps, we obtained that $u(x,t)=0$ in $\mathcal{C}_\theta\times(0,T_{k+1})$, where
$$T_{k+1}=(1-T_k)T_1+T_k\rightarrow1.$$
Thus we proved Theorem \ref{thm}.

\section{Proof of the Carleman Inequality}
In this section, we are aimed to prove the Carleman Inequality (\ref{CI}) in Proposition \ref{Prop}.

First of all, we make the following notations for simplicity.
Denote by $A^T$ the transpose of a matrix $A$, and $x=(x_1, x_2,\dots,x_n)^T$, $e=(1,0,\dots,0)^T$.

Let $\Lambda(t)=t^{-K}-1$, $h=\frac{x_1}{r}$, $f(h)=h^m-\varepsilon^m$, then $\varphi(x)=r^\alpha f(h)$ by the  assumptions of Proposition \ref{Prop}.

We write $$\Phi=a\Lambda(t)\varphi(x)-\frac{|x|^2+K}{16t}=\Phi_1+\Phi_2,$$
where
$$\Phi_1=a\Lambda(t)\varphi(x),\quad \Phi_2=-\frac{|x|^2+K}{16t}.$$
Moreover, $D^2\Phi$, $I_n$ denote the Hessian matrix of $\Phi$, the identity matrix of $\mathbb{R}^n$, respectively.

Under the  assumptions as above and Proposition \ref{Prop}, we derive the following lemma:
\begin{Lemma}\label{lem:bound of f}
For $2<m<3$ and $1<\alpha\leq 2$,
direct calculations show that $f$ has the following properties in the interval $h\in[\varepsilon,1]$:
\begin{equation}\label{p.f}
\begin{split}
&(i)~f(h)\geq 0,~~f(\varepsilon)=0;\\
&(ii)~f''(h)\geq0;\\
&(iii)~(\alpha^2-2\alpha)f(h)+(3-2\alpha)hf'(h)+h^2f''(h)\geq0;\\
&(iv)~(\alpha-1)^2f'(h)^2+(2\alpha-\alpha^2)f(h)f''(h)-hf'(h)f''(h)\leq0.
\end{split}
\end{equation}
\end{Lemma}
{\it Proof of Lemma \ref{lem:bound of f}:}

Since $f(h)=h^m-\varepsilon^m$, we have $f'(h)=mh^{m-1}$ and $f''(h)=(m^2-m)h^{m-2}.$ Obviously $(i)$ and $(ii)$ hold.

For $(iii)$, by $2<m<3$ and $1<\alpha\leq 2$ we have
\begin{equation*}\label{eq:iii}
\begin{split}
&(\alpha^2-2\alpha)f(h)+(3-2\alpha)hf'(h)+h^2f''(h)\\
&=[(m-\alpha)^2+2(m-\alpha)]h^m+(2\alpha-\alpha^2)\varepsilon^m>0
\end{split}
\end{equation*}

Consider $(iv)$, and we get
\begin{equation*}\label{eq:iv}
\begin{split}
&(\alpha-1)^2f'(h)^2+(2\alpha-\alpha^2)f(h)f''(h)-hf'(h)f''(h)\\
&=m[(2m-m^2)-(2\alpha-\alpha^2)]h^{2m-2}-(2\alpha-\alpha^2)(m^2-m)h^{m-2}\varepsilon^m<0
\end{split}
\end{equation*}
The proof is complete.

In the following, we often write $f(h)$, $f'(h)$, $f''(h)$ as $f$, $f'$, $f''$.\\

\emph{\bf Proof of Proposition \ref{Prop}.} Let $u$ be an arbitrary function in $C_0^\infty(Q)$ and $v=e^\Phi u$, then
$$
\int_Qe^{2\Phi}|\partial_t u+\Delta u|^2dxdt
=\int_Q|\partial_t v-2\nabla\Phi\cdot\nabla v+\Delta v+(|\nabla\Phi|^2-\partial_t\Phi-\Delta\Phi)v|^2dxdt.
$$
By the Cauchy inequality we have
$$
\int_Qe^{2\Phi}|\partial_t u+\Delta u|^2dxdt
\geq2\int_Q[\partial_t v-2\nabla\Phi\cdot\nabla v-(\Delta\Phi+\frac{F}{2})v][\Delta v+(|\nabla\Phi|^2-\partial_t\Phi+\frac{F}{2})v]dxdt,
$$
where $F$ is an arbitrary twice differentiable function, to be decided later. Rewrite the right side of the above inequality by integration by parts, then we have
\begin{equation}\label{CIconstant}
\begin{split}
&\int_Qe^{2\Phi}|\partial_t u+\Delta u|^2dxdt\\
\geq&\int_Q[4D^2\Phi\nabla v\cdot\nabla v+F|\nabla v|^2]dxdt\\
&+\int_Q[4D^2\Phi\nabla\Phi\cdot\nabla\Phi-2\partial_t|\nabla\Phi|^2+\partial_{tt}\Phi-\Delta^2\Phi]v^2dxdt\\
&+\int_Q[-\frac{1}{2}\partial_tF-\frac{1}{2}\Delta F+\nabla\Phi\cdot\nabla F-F(|\nabla\Phi|^2-\partial_t\Phi+\frac{F}{2})]v^2dxdt.
\end{split}
\end{equation}
Unfortunately, $D^2\Phi$ is not positive, and we need to choose an appropriate function $F$ to compensate it. Next, firstly we estimate $D^2\Phi$ and choose a suitable function $F$ to make the matrix $4D^2\Phi+FI_n$ be positive; secondly, we estimate the terms including $v^2$.

{\bf Step 1. Estimate for the gradient terms.} In fact, we have the following conclusion.

Claim that: for $F=-4a \Lambda r^{\alpha-2}(\alpha f-hf')+\frac{3}{t}$, there holds
\begin{equation}\label{CI.1-}
\begin{split}
\int_Q[4D^2\Phi\nabla v\cdot\nabla v+F|\nabla v|^2]dxdt
\geq\int_Q\frac{5}{2t}(\frac{1}{2}e^{2\Phi}|\nabla u|^2-|\nabla\Phi|^2v^2)dxdt.
\end{split}
\end{equation}

Recall that $\varphi=r^{\alpha}f$, $\nabla h=r^{-1}e-r^{-2}hx$, then
\begin{equation}\label{ccl1}
\nabla\varphi=\alpha r^{\alpha-2}fx+r^{\alpha}f'\nabla h=r^{\alpha-2}(\alpha f-hf')x+r^{\alpha-1}f'e;
\end{equation}
\begin{equation}\label{ccl2}
\begin{split}
D^2\varphi=&(\alpha-2)r^{\alpha-4}(\alpha f-hf')x x^T+r^{\alpha-2}[(\alpha-1)f'-hf'']x(\nabla h)^T\\
&+r^{\alpha-2}(\alpha f-hf')I_n+(\alpha-1)r^{\alpha-3}f'e x^T+r^{\alpha-1}f''e (\nabla h)^T\\
=&r^{\alpha-2}[(\alpha f-hf')I_n+B],
\end{split}
\end{equation}
where
\begin{equation}\label{ccl3}
\begin{split}
B\equiv&f''e e^T+r^{-1}[(\alpha-1)f'-hf''](e x^T+xe^T)\\
&+r^{-2}[(\alpha^2-2\alpha)f+(3-2\alpha)hf'+h^2f'']xx^T.
\end{split}
\end{equation}
Using $ii)$, $iii)$, $iv)$ of (\ref{p.f}), and
\begin{eqnarray*}
\begin{split}
&r^{-2}[(\alpha-1)f'-hf'']^2-f''r^{-2}[(\alpha^2-2\alpha)f+(3-2\alpha)hf'+h^2f'']\\
=&r^{-2}[(\alpha-1)^2f'^2+(2\alpha-\alpha^2)ff''-hf'f'']\\
\leq&0,
\end{split}
\end{eqnarray*}
we obtain that $B$ is nonnegative by Cauchy inequality. Thus
\begin{equation}\label{hessian.1}
D^2\Phi_1=a\Lambda r^{\alpha-2}[(\alpha f-hf')I_n+B]
\geq a \Lambda r^{\alpha-2}(\alpha f-hf')I_n
\equiv HI_n.
\end{equation}
It is easy to verify that
\begin{equation}\label{p.H}
\begin{split}
&(i)~-C(\varepsilon)a \Lambda r^{\alpha-2}\leq H<0,\\
&(ii)~|\nabla^kH|\leq C(\varepsilon)a \Lambda r^{\alpha-2-k},~~k=1,2,
\end{split}
\end{equation}
and by (\ref{ccl3})
\begin{equation}\label{hessian.2}
-C(\varepsilon)a\Lambda r^{\alpha-2}\leq D^2\Phi_1\leq C(\varepsilon)a \Lambda r^{\alpha-2}.
\end{equation}
Immediately from (\ref{hessian.1}) we get
\begin{equation}\label{hessian.3}
D^2\Phi=D^2\Phi_1-\frac{1}{8t}I_n\geq (H-\frac{1}{8t})I_n.
\end{equation}
We choose
\begin{equation}\label{F}
F=-4H+\frac{3}{t},
\end{equation}
then by (\ref{hessian.3}), (\ref{F}) and Cauchy inequality, we obtain that
\begin{equation*}\label{CI.1}
\begin{split}
\int_Q[4D^2\Phi\nabla v\cdot\nabla v+F|\nabla v|^2]dxdt
\geq&\int_Q\frac{5}{2t}|\nabla v|^2dxdt\\
=&\int_Q\frac{5}{2t}|e^\Phi\nabla u+\nabla\Phi v|^2dxdt\\
\geq&\int_Q\frac{5}{2t}(\frac{1}{2}e^{2\Phi}|\nabla u|^2-|\nabla\Phi|^2v^2)dxdt.
\end{split}
\end{equation*}
Hence, the proof of (\ref{CI.1-}) is complete.

{\bf Step 2. Estimate for the $v^2$ terms.}
By (\ref{CIconstant}) and (\ref{CI.1-}), we obtain
\begin{equation}\label{CI.2}
\int_Qe^{2\Phi}|\partial_t u+\Delta u|^2dxdt\geq\int_Q\frac{5}{4t}e^{2\Phi}|\nabla u|^2dxdt+\int_QJv^2dxdt,
\end{equation}
where
\begin{eqnarray*}
\begin{split}
J=&4D^2\Phi\nabla\Phi\cdot\nabla\Phi-2\partial_t|\nabla\Phi|^2+\partial_{tt}\Phi-\Delta^2\Phi-\frac{5}{2t}|\nabla\Phi|^2\\
&-\frac{1}{2}\partial_tF-\frac{1}{2}\Delta F+\nabla\Phi\cdot\nabla F-F(|\nabla\Phi|^2-\partial_t\Phi+\frac{F}{2}).
\end{split}
\end{eqnarray*}

In order to estimate $J$, we divide $J$ into four parts according to the orders of the parameter $a$ and will estimate  each part of them later. Using the definition  of $\Phi_1$ and (\ref{ccl1}),  we have
$$|\nabla\Phi|^2=|\nabla\Phi_1|^2+|\nabla\Phi_2|^2+2\nabla\Phi_1\cdot\nabla\Phi_2=|\nabla\Phi_1|^2+|\nabla\Phi_2|^2-\frac{\alpha}{4t}\Phi_1,$$
and
$D^2\Phi_2\nabla\Phi\cdot\nabla\Phi=-\frac{1}{8t}|\nabla\Phi|^2,$
then by (\ref{F}) we get
$$J=J_3+J_2+J_1+J_0,$$
where
\begin{equation}\label{J}
\begin{split}
J_3=&4D^2\Phi_1\nabla\Phi_1\cdot\nabla\Phi_1+4H|\nabla\Phi_1|^2;\\
J_2=&8D^2\Phi_1\nabla\Phi_1\cdot\nabla\Phi_2-2\partial_t|\nabla\Phi_1|^2-4H\partial_t\Phi_1-\frac{6}{t}|\nabla\Phi_1|^2-4\nabla\Phi_1\cdot\nabla H\\
&-\frac{\alpha}{t}H\Phi_1-8H^2;\\
J_1=&4D^2\Phi_1\nabla\Phi_2\cdot\nabla\Phi_2+4H|\nabla\Phi_2|^2-4\nabla\Phi_2\cdot\nabla H+\frac{12}{t}H+2\Delta H-\Delta^2\Phi_1\\
&+2\partial_tH-4H\partial_t\Phi_2+\partial_{tt}\Phi_1+\frac{\alpha}{2}\partial_t(\frac{\Phi_1}{t})+\frac{3\alpha}{2t^2}\Phi_1+\frac{3}{t}\partial_t\Phi_1;\\
J_0=&\partial_{tt}\Phi_2+\frac{3}{t}\partial_t\Phi_2-2\partial_t|\nabla\Phi_2|^2-\frac{6}{t}|\nabla\Phi_2|^2-\frac{3}{t^2}.
\end{split}
\end{equation}

\emph{\bf Estimate of $J_3$.}\\

Due to (\ref{ccl1})-(\ref{ccl3}), the definition of $\Phi_1$ and $H$, we have
\begin{eqnarray*}
J_3=4a^3\Lambda^3D^2\varphi\nabla\varphi\cdot\nabla\varphi+4a^3\Lambda^3r^{\alpha-2}(\alpha f-hf')|\nabla\varphi|^2
\equiv4a^3\Lambda^3r^{3\alpha-4}l_1(\alpha,h),
\end{eqnarray*}
where
\begin{eqnarray*}
l_1(\alpha,h)=\alpha^4f^3-\alpha^2hf^2f'+2\alpha^2(1-h^2)ff'^2
-2h(1-h^2)f'^3+(1-h^2)^2f'^2f''.
\end{eqnarray*}
We hope that
\begin{equation}\label{l_1alpha-h}
l_1(\alpha,h)\geq0.
\end{equation}
For $\frac12<\gamma\leq 1$, we estimate $l_1(\alpha,h)$ as follows:
\begin{eqnarray*}
\begin{split}
l_1(\alpha,h)=&f(\alpha^2f-\gamma hf')^2+hf f'((2\gamma-1)\alpha^2f-\gamma^2hf')\\
&+2(1-h^2)f'^2(\alpha^2f-hf'+\frac{1-h^2}{2}f'')\\
\geq&hf f'((2\gamma-1)\alpha^2f-\gamma^2hf')+2(1-h^2)f'^2(\alpha^2f-hf'+\frac{1-h^2}{2}f'')\\
=&hf f'[(2\gamma-1)\alpha^2f-\gamma^2hf'+\gamma^2\frac{m-1}{2h}(1-h^2)f']\\
&+2(1-h^2)f'^2[(\alpha^2-\gamma^2\frac{m-1}{4})f-hf'+\frac{1-h^2}{2}f'']
\end{split}
\end{eqnarray*}
Since $f(h)=h^m-\varepsilon^m$, we have $f''=\frac{m-1}{h}f'$, then
\begin{eqnarray*}
\begin{split}
l_1(\alpha,h)\geq &hf f'[(2\gamma-1)\alpha^2f-\gamma^2hf'+\gamma^2\frac{1-h^2}{2}f'']\\
&+2(1-h^2)f'^2[(\alpha^2-\gamma^2\frac{m-1}{4})f-hf'+\frac{1-h^2}{2}f'']
\end{split}
\end{eqnarray*}
Hence the condition $(\ref{l_1alpha-h})$ is satisfied if the following two inequalities hold:
\begin{equation}\label{eq:2gamma-1}
(2\gamma-1)\alpha^2\geq \gamma^2(\alpha^2-\gamma^2\frac{m-1}{4}),
\end{equation}
and for $h\in [\varepsilon,1]$,
\begin{eqnarray}\label{l2alpha-h}
\begin{split}
l_2(\alpha,h)\equiv&(\alpha^2-\gamma^2\frac{m-1}{4})f-hf'+\frac{1-h^2}{2}f''\\
=&(\alpha^2-\gamma^2\frac{m-1}{4}-\frac{m^2+m}{2})h^m+\frac{m^2-m}{2}h^{m-2}-(\alpha^2-\gamma^2\frac{m-1}{4})\varepsilon^m\geq 0.
\end{split}
\end{eqnarray}
Notice that for $\gamma>\frac12$, $1<\alpha\leq 2$ and $2.36\leq m<3$, we have
\begin{equation*}
\alpha^2-\gamma^2\frac{m-1}{4}-\frac{m^2+m}{2}\leq 0,
\end{equation*}
consequently,
\begin{equation}\label{j1.1}
\frac{d^2}{dh^2}l_2(\alpha,h)<0.
\end{equation}
Hence $l_2(\alpha,h)$ is a concave function, and $l_2(\alpha,\varepsilon), l_2(\alpha,1)\geq 0$ yield that $(\ref{l2alpha-h})$. Moreover, for $\varepsilon\in(0,\varepsilon_0)$,
\begin{equation}\label{j1.2}
l_2(\alpha,\varepsilon)=\frac{m(m+1)}{2}\varepsilon^{m-2}(\frac{m-1}{m+1}-\varepsilon^2)
>\frac{m(m+1)}{2}\varepsilon^{m-2}(\frac{m-1}{m+1}-\varepsilon_0^2)
\geq0,
\end{equation}
$$l_2(\alpha,1)=(\alpha^2-\gamma^2\frac{m-1}{4}-m)-(\alpha^2-\gamma^2\frac{m-1}{4})\varepsilon^m> 0.$$
Since $\varepsilon$ is a increasing function with respect to $\alpha,$  to obtain a bigger $\varepsilon_0$ we only need to consider $\alpha\rightarrow2.$
Concluding the above arguments $(\ref{eq:2gamma-1})$ and $(\ref{j1.2})$, in order to ensure that the result $(\ref{l_1alpha-h})$ holds, we are going to seek the following sharpest $\varepsilon_0$ such that for $\frac12<\gamma\leq 1$, $2.36\leq m<3$, there hold
\begin{eqnarray}\label{eq:j3 condition}
\left\{\begin{array}{lll}
&(i)~4(2\gamma-1)\geq \gamma^2(4-\gamma^2\frac{m-1}{4});\\
&(ii)~\frac{m-1}{m+1}-\varepsilon_0^2\geq 0;\\
&(iii)~(4-\gamma^2\frac{m-1}{4}-m)-(4-\gamma^2\frac{m-1}{4})\varepsilon_0^m\geq 0.
\end{array}\right.
\end{eqnarray}
Especially, when $\gamma=1$, $(i)$ holds obviously and $(ii)$-$(iii)$ show that $(m,\varepsilon_0)$ the intersection point of the two functions $g_1(p)=\sqrt{\frac{p-1}{p+1}}$ and $g_2(p)=(\frac{17-5p}{17-p})^{\frac{1}{p}}$ in the interval $p\in[2.36,3]$. Then $m\approx2.39$, $\varepsilon_0=\sqrt{\frac{m-1}{m+1}}\approx 0.64$.

To obtain a sharper $\varepsilon_0$, we consider $\gamma<1,$ and we could solve the above three inequalities by Newton approximate methods.
Indeed, by (\ref{eq:j3 condition}) we derive that $m\approx 2.4600$, $\varepsilon_0=\sqrt{\frac{m-1}{m+1}}\approx 0.6495$ and $\gamma\approx0.8092$. Hence $\theta_0=\frac{\arccos(0.6495)\times 360}{3.1415927}\approx 98.99^{\circ}.$

Finally, for $\varepsilon<\varepsilon_0$
there exists $\delta_1(\varepsilon)<1$, such that when $\alpha\in (2-\delta_1,2)$, $(\ref{eq:2gamma-1})$ and $(\ref{j1.2})$ hold. Hence
\begin{equation}\label{j1.3}
l_1(\alpha,h)\geq 0.
\end{equation}
Then, for any $\varepsilon\in(0,\varepsilon_0)$, choosing $\alpha<2$ suitably  we have
\begin{equation}\label{J_3}
J_3\geq0.
\end{equation}

\emph{\bf Estimate of $J_2$.}\\

By (\ref{hessian.2}) we know that
$$D^2\Phi_1\nabla\Phi_1\cdot\nabla\Phi_2\geq-C(\varepsilon)a \Lambda r^{\alpha-2}|\nabla \Phi_1||\nabla \Phi_2|;$$ and by (\ref{p.H}) $H<0$, thus $-\frac{\alpha}{t}H\Phi_1\geq 0$. Then by (\ref{J}) we have
\begin{equation}\label{J_2.1}
\begin{split}
J_2\geq& -C(\varepsilon)a \Lambda r^{\alpha-2}|\nabla \Phi_1||\nabla \Phi_2|-2\partial_t|\nabla\Phi_1|^2-4H\partial_t\Phi_1\\
&-\frac{6}{t}|\nabla\Phi_1|^2-4|\nabla\Phi_1||\nabla H|-8H^2.
\end{split}
\end{equation}
Next we estimate the terms on the right side of the above inequality. \\
First, (\ref{ccl1}) and direct calculations yield that
\begin{equation}\label{estofpsi1}
|\nabla\Phi_1|=a \Lambda|\nabla\varphi|=a \Lambda r^{\alpha-1}\sqrt{\alpha^2f^2+(1-h^2)f'^2},
\end{equation}
and it is easy to verify that
\begin{equation}\label{estofpsi2}
|\nabla\Phi_1|\leq C(\varepsilon)a\Lambda r^{\alpha-1}.
\end{equation}
By (\ref{p.H}) and (\ref{estofpsi2}), we have
\begin{equation}\label{J_2.2}
\begin{split}
-C(\varepsilon)a\Lambda r^{\alpha-2}|\nabla \Phi_1||\nabla \Phi_2|\geq-\frac{C(\varepsilon)}{t}a^2\Lambda^2 r^{2\alpha-2};\\
-\frac{6}{t}|\nabla\Phi_1|^2\geq-\frac{C(\varepsilon)}{t}a^2\Lambda^2 r^{2\alpha-2};\\
-4|\nabla\Phi_1||\nabla H|\geq-C(\varepsilon)a^2\Lambda^2 r^{2\alpha-4};\\
-8H^2\geq-C(\varepsilon)a^2\Lambda^2 r^{2\alpha-4}.
\end{split}
\end{equation}
Recall  that $H=a\Lambda r^{\alpha-2}(\alpha f-hf')$ by (\ref{hessian.1}), and from (\ref{estofpsi1}) we derive that
\begin{equation}\label{J_2.3}
\begin{split}
&-2\partial_t|\nabla\Phi_1|^2-4H\partial_t\Phi_1\\
=&-4a^2\Lambda\Lambda'r^{2\alpha-2}[\alpha^2f^2+(1-h^2)f'^2]-4a^2 \Lambda\Lambda'r^{2\alpha-2}f(\alpha f-hf')\\
=&-4a^2\Lambda\Lambda'r^{2\alpha-2}[\alpha^2f^2+(1-h^2)f'^2+f(\alpha f-hf')],
\end{split}
\end{equation}
then combining (\ref{J_2.1}), (\ref{J_2.2}) and (\ref{J_2.3}) we have (note that $t\in(0,1)$ and $r>1$)
\begin{equation}\label{estJ_2}
\begin{split}
J_2\geq&-\frac{C(\varepsilon)}{t}a^2\Lambda^2 r^{2\alpha-2}-4a^2\Lambda\Lambda'r^{2\alpha-2}[(\alpha^2+\alpha)f^2-hff'+(1-h^2)f'^2]\\
\equiv&-\frac{C(\varepsilon)}{t}a^2\Lambda^2 r^{2\alpha-2}-4a^2\Lambda\Lambda'r^{2\alpha-2}l_3(\alpha,h),
\end{split}
\end{equation}
and
\begin{eqnarray*}
\begin{split}
l_3(\alpha,h)=&(\alpha^2+\alpha)f^2-hff'+(1-h^2)f'^2\\
=&(\alpha^2+\alpha)\varepsilon^{2m}+h^m[(\alpha^2+\alpha-m^2-m)h^m+m^2h^{m-2}-(2\alpha^2+2\alpha-m)\varepsilon^m]\\
\equiv&(\alpha^2+\alpha)\varepsilon^{2m}+h^ml_4(\alpha,h),
\end{split}
\end{eqnarray*}
where $$l_4(\alpha,h)=(\alpha^2+\alpha-m^2-m)h^m+m^2h^{m-2}-(2\alpha^2+2\alpha-m)\varepsilon^m.$$
We remark that $(\alpha^2+\alpha-m^2-m)<0$, and $l_4(\alpha,h)$ and $l_2(\alpha,h)$ are of the same type concave function. Similar arguments as $l_2(\alpha,h)$, one prove that for $m=2.46$, $\varepsilon_0=0.6495$ and $\varepsilon\in(0,\varepsilon_0)$, there hold
$$
l_4(2,\varepsilon)>0,\quad l_4(2,1)>0,
$$
then there exists some $\delta_2(\varepsilon)<\delta_1$, such that when $\alpha\in (2-\delta_2,2)$, we have
$$l_4(\alpha,h)\geq0.$$
Hence
\begin{equation}\label{estl_3}
l_3(\alpha,h)\geq(\alpha^2+\alpha)\varepsilon^{2m}\geq2\varepsilon^{2m}.
\end{equation}
Note that $\Lambda'<0$, and by (\ref{estJ_2}) and (\ref{estl_3}) we obtain
\begin{equation}\label{J_2}
J_2\geq-\frac{C(\varepsilon)}{t}a^2\Lambda^2 r^{2\alpha-2}-8\varepsilon^{2m}a^2\Lambda\Lambda'r^{2\alpha-2}.
\end{equation}

\emph{\bf Estimate of $J_1$.}\\

Since $H<0$, $\partial_tH=\frac{\Lambda'}{\Lambda}H\geq0$, $\partial_t\Phi_2>0$ and $D^2\Phi_1\geq HI_n$ by (\ref{hessian.1}), then from (\ref{J}) we obtain that
\begin{eqnarray*}
\begin{split}
J_1\geq&8H|\nabla\Phi_2|^2-4|\nabla\Phi_2||\nabla H|+\frac{12}{t}H+2\Delta H-\Delta^2\Phi_1\\
&+\partial_{tt}\Phi_1+\frac{\alpha}{2}\partial_t(\frac{\Phi_1}{t})+\frac{3\alpha}{2t^2}\Phi_1+\frac{3}{t}\partial_t\Phi_1.
\end{split}
\end{eqnarray*}
From (\ref{p.H}), $|\nabla\Phi_2|\leq \frac{|x|}{8t}$ and
$$\Delta^2\Phi_1=a\Lambda\Delta^2\varphi\leq C(\varepsilon)a \Lambda r^{\alpha-4},$$
we derive that
\begin{equation*}
\begin{split}
J_1\geq&-\frac{C(\varepsilon)}{t^2}a \Lambda r^{\alpha}-\frac{C(\varepsilon)}{t}a\Lambda r^{\alpha-2}-C(\varepsilon)a\Lambda r^{\alpha-4}\\
&+a\varphi[\Lambda''+\frac{\alpha}{2}(\frac{\Lambda'}{t}-\frac{\Lambda}{t^2})+\frac{3\alpha}{2t^2}\Lambda+\frac{3}{t}\Lambda']\\
\geq&-\frac{C(\varepsilon)}{t^2}a\Lambda r^{\alpha}+a\varphi(\Lambda''+\frac{\alpha+6}{2t}\Lambda'+\frac{\alpha}{t^2}\Lambda)\\
\geq&-\frac{C(\varepsilon)}{t^2}a\Lambda r^{\alpha}++a\varphi(\Lambda''+\frac{4}{t}\Lambda'),
\end{split}
\end{equation*}
where we have used $r>1$, $t\in (0,1)$, $\Lambda>0$ and $\Lambda'<0.$

At last, apply Cauchy inequality, and we get
\begin{equation}\label{J_1}
J_1\geq-\frac{C(\varepsilon)}{t}a^2\Lambda^2 r^{2\alpha-2}-\frac{r^2}{32t^3}+a\varphi(\Lambda''+\frac{4}{t}\Lambda').
\end{equation}

\emph{\bf Estimate of $J_0$.}\\

By (\ref{J}), direct calculations show that
\begin{equation}\label{J_0}
\begin{split}
J_0=&-\frac{r^2+K}{8t^3}+\frac{3(r^2+K)}{16t^3}+\frac{r^2}{16t^3}-\frac{3r^2}{32t^3}-\frac{3}{t^2}\\
=&\frac{r^2}{32t^3}+\frac{K}{16t^3}-\frac{3}{t^2}.
\end{split}
\end{equation}
Combing (\ref{J_3}), (\ref{J_2}), (\ref{J_1}) and (\ref{J_0}), we obtain that when $\alpha\in(2-\delta_2,2)$,
\begin{eqnarray*}
\begin{split}
J\geq&-8\varepsilon^{2m}a^2\Lambda\Lambda'r^{2\alpha-2}-\frac{C(\varepsilon)}{t}a^2\Lambda^2 r^{2\alpha-2}
+a\varphi(\Lambda''+\frac{4}{t}\Lambda')+\frac{K}{16t^3}-\frac{3}{t^2}\\
\geq&8a^2\Lambda r^{2\alpha-2}[-\varepsilon^{2m}\Lambda'-\frac{C(\varepsilon)}{t}\Lambda]+a\varphi K(K-3)t^{-K-2}+\frac{K-48}{16t^3}\\
\geq&8a^2\Lambda r^{2\alpha-2}t^{-K-1}[\varepsilon^{2m}K-C(\varepsilon)]+a\varphi K(K-3)t^{-K-2}+\frac{K-48}{16t^3},
\end{split}
\end{eqnarray*}
and choosing $K=K(\varepsilon)$ large enough, we have
$$J\geq\frac{1}{t^3}.$$
Then
\begin{equation}\label{J.over}
\int_QJv^2dxdt\geq \int_Qe^{2\Phi}\frac{u^2}{t^3}dxdt.
\end{equation}
 Combining (\ref{CI.2}) and (\ref{J.over}), finally we obtain that for $m=2.46$, $\varepsilon_0=0.6495$, and $0<\varepsilon<\varepsilon_0$, there exist $\alpha=\alpha(\varepsilon)<2$, and  $K=K(\varepsilon)$ large enough such that
$$\int_Qe^{2\Phi}|\partial_t u+\Delta u|^2dxdt\geq\int_Qe^{2\Phi}(\frac{|\nabla u|^2}{t}+\frac{u^2}{t^3})dxdt.$$
Thus we have completed the proof of  Carleman inequality (\ref{CI}).

\section{Conclusion}

In the arguments of applying Carleman inequality to prove (BU), the construction of weight function $\Phi$ is crucial. As said in \cite{LlS} by Lu Li and V. \v{S}ver\'{a}k that, the main difficulty of (BU) on cones lies in that

\emph{We can never find a positive function $\varphi$ defined on a cone, such that $\varphi=0$ on the boundary of the cone and the Hessian matrix of $\varphi$ is nonnegative. Actually, such kind of function doesn't exist. }

In this paper, we construct $\Phi$ to be  of this type
$$\Phi=\gamma(t^{-K}-1)\varphi(x)-\frac{|x|^2+K}{16t}.$$
Generally speaking, we require that $\varphi=0$ on the boundary of the domain.
Note that the boundary of the cone $\{x|x_1>\varepsilon r\}$ is exactly $\frac{x_1}{r}=\varepsilon$, then we choose
$$\varphi(x)=r^\alpha f(h),~~h=\frac{x_1}{r},~~1<\alpha<2,$$
with $f$ satisfying (\ref{p.f}) first, that is
\begin{eqnarray*}
\begin{split}
\left\{
  \begin{array}{ll}
    f\geq 0,~~f(\varepsilon)=0,  \\
    f''\geq0,  \\
    (\alpha^2-2\alpha)f+(3-2\alpha)hf'+h^2f''\geq0,  \\
    (\alpha-1)^2f'^2+(2\alpha-\alpha^2)ff''-hf'f''\leq0.
  \end{array}
\right.
\end{split}
\end{eqnarray*}
Then the range of $\varepsilon$ and $\theta$ are determined mainly by (\ref{l_1alpha-h}), that is
$$\alpha^4f^3-\alpha^2hf^2f'+2\alpha^2(1-h^2)ff'^2-2h(1-h^2)f'^3+(1-h^2)^2f'^2f''\geq0.$$

If we choose $f$ to be of the simplest type
$$f(h)=h^\beta-\varepsilon^\beta,$$
then it is easy to verify that (\ref{p.f}) holds when $\beta\geq\alpha$.
In particular, when $\beta=\alpha$, it responds to the weight function constructed by Lu Li and V. \v{S}ver\'{a}k in \cite{LlS} and in this case
$$\varepsilon\in(0,\sqrt{\frac{1}{3}}),~~\theta\in(2\arccos\sqrt{\frac{1}{3}},\pi);$$
 when $\beta=m>\alpha$, it responds to the weight function constructed in this paper,
$$\varepsilon\in(0,\varepsilon_0),~~\theta\in(2\arccos{\varepsilon_0},\pi),~~\varepsilon_0\approx0.6495, $$
and we obtain a wider range of $\varepsilon$ and $\theta$.

This means that: if we choose finer and more complicated function $f$, which satisfies (\ref{p.f}) at least, then the range of $\varepsilon$ and $\theta$ which is determined by (\ref{l_1alpha-h}) will become correspondingly wider.

However, can we find a function $f$, such that the range of $\varepsilon$ and $\theta$ we obtained are exactly
$$\varepsilon\in(0,\sqrt{\frac{1}{2}}),~~\theta\in(\frac{\pi}{2},\pi)?$$
It seems to be somewhat difficult because of the limitations of the methods and the choosing of most suitable function which one has no idea. If the (BU) problem on cones could be solved completely, some new ideas or observations are necessary.

\bigskip

\noindent {\bf Acknowledgments.}
The authors would like to thank Professor Liqun Zhang for many valuable discussions on this topic. Wang is supported by "the Fundamental Research Funds for the Central Universities" and partially by The Institute of Mathematical Sciences of CUHK.

\vspace{1cm}

{\small}

\noindent

{\bf Jie Wu}\\
Beijing International Center for Mathematical Research, Peking University, 5 Yiheyuan Road, Haidian District, Beijing 100871, China. \\
Email address: {\bf jackwu@amss.ac.cn}\\

{\bf Wengdong Wang}\\
School of  Mathematical Sciences, Dalian University of Technology, Dalian 116024, P.R. China. \\
Email address: {\bf wendong@dlut.edu.cn}

\end{document}